\documentclass[11pt]{amsart}
\usepackage{amssymb, latexsym}

\theoremstyle{plain}
\newtheorem{theorem}{Theorem}
\newtheorem{corollary}{Corollary}

\newtheorem {lemma}{Lemma}
\newtheorem{proposition}{Proposition}

\theoremstyle{remark}

\newtheorem*{Remark 1}{Remark 1}
\newtheorem*{Remark 2}{Remark 2}
\newtheorem*{Remark 3}{Remark 3}
\newtheorem*{Remark 4}{Remark 4}

\numberwithin{equation}{section}

\begin{document}

\title[Detecting Tampering]%
 { Detecting Tampering in a Random Hypercube}

\author{ Ross G. Pinsky}
\address{Department of Mathematics\\
Technion---Israel Institute of Technology\\
Haifa, 32000\\ Israel} \email{ pinsky@math.technion.ac.il}
\urladdr{http://www.math.technion.ac.il/~pinsky/}

\date{}

\begin{abstract}
Consider the random hypercube $H_2^n(p_n)$ obtained from the hypercube $H_2^n$ by deleting any given edge with probabilty
$1-p_n$, independently of all the other edges. A diameter path in $H_2^n$ is a longest geodesic path   in $H_2^n$.
Consider the following two ways of tampering with the random graph $H_2^n(p_n)$: (i) choose a diameter path at random
and adjoin all of its edges to $H_2^n(p_n)$; (ii)  choose a diameter path at random
from among those that start at $0=(0,\cdots, 0)$, and adjoin all of its edges to $H_2^n(p_n)$.
We study the question of  whether these tamperings are detectable asymptotically as $n\to\infty$.
\end{abstract}

\maketitle
\section{Introduction and Statement of Results}\label{S:intro}

Let $H_2^n=(V_n,e_n)$ denote the
 $n$-dimensional hypercube.
 Recall that  the vertices $V_n$ of $H_2^n$ are identified with   $\{0,1\}^n$,
and  an edge in $e_n$ connects two vertices if and only if they differ in exactly one component.
Denote vertices by $\bar x=(x_1,\cdots,x_n)$.
A geodesic path from $\bar x$ to $\bar y$ is a shortest path from $\bar x$ to $\bar y$. A diameter path in $H_2^n$ is a longest geodesic path
in $H_2^n$. The set of diameter paths is the set of paths $\bar x_0\bar x_1\cdots \bar x_n$, where
$\bar x_n=\bar 1-\bar x_0$ and  $\bar 1\equiv (1,\cdots, 1)$.

Let $H_2^n(p_n)$  denote the  random hypercube obtained by starting with  the graph  $H_2^n$ and deleting
any given edge with probability  $1-p_n$, independently of all the other edges.
Let $P_{n,p_n}$ denote the corresponding probability measure;  $P_{n,p_n}$ is a measure on
$\mathcal{E}_n\equiv2^{e_n}$, the space of all subsets of $e_n$. An  element of $\mathcal{E}_n$
will be called an edge configuration.

We consider two similar ways of  tampering with the random hypercube. The first way is to
 choose
 a diameter path from $H_2^n$ at random and adjoin it to $H_2^n(p_n)$; that is, we ``add'' to the random graph every edge of  this
diameter path that is  not already in  the random graph.
Denote the induced measure on $\mathcal{E}_n$ by $P_{n,p_n}^{\text{tam}}$.
The second way is to consider $0\equiv(0,\cdots, 0)$ as a distinguished vertex in the  hypercube, and to adjoin to
the random hypercube a diameter path chosen at random from among
those diameter paths which start at 0. Denote the induced measure on $\mathcal{E}_n$ by $P_{n,p_n}^{\text{tam},0}$.

Can one detect the tampering asymptotically as $n\to\infty$?  Let $Q_n$ be generic notation for either
$P_{n,p_n}^{\text{tam}}$ or $P_{n,p_n}^{\text{tam},0}$.
Let $||P_{n,p_n}-Q_n||_{\text{TV}}$ denote the total variation distance between the probability measures $P_{n,p_n}$ and $Q_n$.
If $\lim_{n\to\infty}||P_{n,p_n}-Q_n||_{\text{TV}}=1$,  we call the tampering   \it detectable\rm.
If \linebreak $\lim_{n\to\infty}||P_{n,p_n}-Q_n||_{\text{TV}}=0$,  we call the tampering \it  strongly undetectable\rm, while
if $\{||P_{n,p_n}-Q_n||_{\text{TV}}\}_{n=1}^\infty$ is bounded away  from 0 and 1,  we  call the tampering  \it weakly undetectable\rm.

The number of diameter paths in $H_2^n$ is easily seen to be $2^{n-1}n!$, while the number of diameter paths
in $H_2^n$ that start from 0 is $n!$. Let $m_n$ denote the number of diameter paths in either of these
two cases. Numbering the diameter paths from 1 to $m_n$, let $O_{n,j}$ denote  the set of edge configurations which contain the $j$-th diameter path.
From the above description of the tampered measures $Q_n=P_{n,p_n}^{\text{tam}}$ or $Q_n=P_{n,p_n}^{\text{tam},0}$, it follows that
\begin{equation}\label{tamperedmeas}
Q_n(\cdot)\equiv \frac1{m_n}\sum_{j=1}^{m_n} P_{n,p_n}(\cdot \ |O_{n,j}).
\end{equation}

Let $N_n^{\text{diam}}:\mathcal{E}_n\to \{0,1,\cdots,m_n\}$ denote the number of diameter paths
in an edge configuration, and let
 $N_n^{\text{diam},0}:\mathcal{E}_n\to  \{0,1,\cdots,m_n\}$ denote the number of diameter paths starting from 0
in an edge configuration. Let $N_n$ be generic notation for either
 $N_n^{\text{diam}}$ or $N_n^{\text{diam},0}$. We refer to $N_n$ as the diameter counting function.

The following proposition, which we prove in the next section, shows that the tampered measure is in fact obtained from the original measure
by size biasing with respect to the diameter counting function $N_n$.

\begin{proposition}\label{sizebias}
Let $Q_n$ denote either of the two tampered measures, and let $N_n$ denote the corresponding  diameter counting function. Then
$$
Q_n(\omega)=\frac{N_n(\omega)}{E_{n,p_n}N_n}P_{n,p_n}(\omega), \ \omega\in \mathcal{E}_n.
$$
\end{proposition}
The following proposition is immediate in light of Proposition \ref{sizebias}.
\begin{proposition}\label{2ndmom}
Let $Q_n$ denote either of the two tampered measures, and let $N_n$ denote the corresponding  diameter counting function. Then
$$
\lim_{n\to\infty}||P_{n, p_n}-Q_n||_{\text{TV}}=0
$$
if and only if the weak law of numbers holds for $N_n$ under $P_{n,p_n}$; that is, if and only if
$$
\lim_{n\to\infty}P_{n,p_n}(|\frac{N_n}{E_{n,p_n}N_n}-1|>\epsilon)=0,\  \text{for all}\ \epsilon>0.
$$
\end{proposition}
The second moment method then yields  the following corollary. Let $\text{Var}_{n,p_n}$ denote the variance with respect to $P_{n,p_n}$.
\begin{corollary}\label{cheb}
Let $Q_n$ denote either of the two tampered measures, and let $N_n$ denote the corresponding  diameter counting function.

\noindent i. If $\text{Var}_{n,p_n}(N_n)=o\big((E_{n,p_n}N_n)^2\big)$,
then $\lim_{n\to\infty}||P_{n,p_n}-Q_n||_{\text{TV}}=0$ and the tampering is strongly undetectable;

\noindent ii.  If $\text{Var}_{n,p_n}(N_n)=O\big((E_{n,p_n}N_n)^2\big)$,
then $\{||P_{n,p_n}-Q_n||_{\text{TV}}\}_{n=1}^\infty$ is bounded away from 1; thus the tampering is not detectable.
\end{corollary}
Part (i) of the corollary  of course follows from Chebyshev's inequality;  we give a  proof of part (ii) in section 2.

We will prove the following result.
\begin{theorem}\label{hyper}
\noindent a. Consider the random hypercube $H_2^n(p_n)$ and tamper with it by adding a random  diameter path. Let
$N_n^{\text{diam}}$ denote the diameter counting function.

\noindent i. If $p_n\le \frac \gamma n$, with $\gamma<\frac e2$, then
 the tampering is detectable; furthermore, the distribution of $N_n^{\text{diam}}$
under $P_{n,p_n}$ converges to the $\delta$-distribution at 0;

\noindent ii. If $p_n\ge \frac \gamma n$, with $\gamma>\frac e2$, then
 the tampering is strongly undetectable; equivalently, the distribution of $N_n^{\text{diam}}$
under $P_{n,p_n}$ satisfies the law of large numbers.

\noindent b. Consider the random hypercube $H_2^n(p_n)$ and tamper with it by adding a random  diameter path that starts from 0.
Let $N_n^{\text{diam},0}$ denote the diameter counting function.

\noindent i. If $p_n\le \frac \gamma n$, with $\gamma<e$, then
 the tampering is detectable; furthermore, the distribution of $N_n^{\text{diam},0}$
under $P_{n,p_n}$ converges to the $\delta$-distribution at 0;

\noindent ii.  If $p_n\ge \frac \gamma n$, with $\gamma> e$, and $\limsup_{n\to\infty}np_n<\infty$, then the tampering is weakly undetectable;
in particular, the distribution of $N_n^{\text{diam},0}$
under $P_{n,p_n}$ does not satisfy  the law of large numbers;

\noindent iii. If $\lim_{n\to\infty}np_n=\infty$,
then the tampering is strongly undetectable; equivalently, the distribution of $N_n^{\text{diam},0}$
under $P_{n,p_n}$ satisfies the law of large numbers.
\end{theorem}

\bf\noindent Remark.\rm\ If  under $P_{n,p_n}$, the distribution of $N_n$
converges to the $\delta$-distribution at 0, then the tampering is detectable
since under the tampered measure one has $N_n\ge1$ a.s.
 By Proposition \ref{sizebias}, if the tampering is strongly
undetectable, then the distribution of $N_n$ must converge to the $\delta$-function at $\infty$.
Naive intuition might suggest that for a tampering problem of the above type, the above two statements
should be if and only if statements, except perhaps conceivably in some narrow bifurcation region between two regimes.
 Theorem \ref{hyper} shows that this is indeed the case for the tampering problem under consideration.
(The proof of the theorem will reveal that in case  (b-ii), the distribution of $N_n^{\text{diam},0}$
converges neither to the $\delta$-distribution at 0 nor to the $\delta$-distribution at $\infty$.)
However, we now point out two examples of similar tampering problems where this intuition fails.
\medskip

\noindent \it Example 1.\rm\ Let $G(n)$ be   the complete graph on $n$ vertices, and let
$G(n,p_n)$ be the Erdos-Renyi random graph with edge probabilities $p_n$; that is, $G(n,p_n)$ is obtained from $G(n)$ by
deleting any particular edge with probability $1-p_n$, independently of all the other edges.
Let $P_{n,p_n}$ denote the corresponding probability measure on edge configurations. As above,
denote the space of all edges by $e_n$ and the space of  all possible edge configurations by $\mathcal{E}_n$.
 Recall that a\it\ Hamiltonian path\rm\ in $G(n)$
 is a path that traverses each of the vertices of the graph exactly once; that is, a path of the form
$x_1x_2\cdots x_n$, where  the $x_i$ are all distinct.
Tamper with the random graph by choosing at random
 a Hamiltonian path from $G(n)$ and adjoining  it to $G(n,p_n)$; that is, ``add'' to the random graph every edge of  this
Hamiltonian  path that is  not already in  the random graph. Call the induced measure $P_{n,p_n}^\text{Ham}$.
 The number of Hamiltonian paths in $G(n)$ is $m_n\equiv\frac12n!$.
 Let $N_n^{\text{ham}}:\mathcal{E}_n\to \{0,1,\cdots, m_n\}$ denote the number of Hamiltonian paths in an edge configuration; we call $N_n^{\text{ham}}$ the
Hamiltonian path counting function.
Quite sophisticated graph theoretical techniques along with probabilistic analysis have yielded  the following beautiful result:
if $p_n=\frac{\log n+\log\log n+\omega_n}n$,  then
\begin{equation}\label{ham}
\begin{aligned}
&\lim_{n\to\infty}P_{n,p_n}(N_n^{\text{Ham}}\ge k)=1, \ \text{for all}\ k, \ \text{if}\ \lim_{n\to\infty}\omega_n=\infty;\\
&\lim_{n\to\infty}P_{n,p_n}(N_n^{\text{Ham}}=0)=1, \ \text{if}\ \lim_{n\to\infty}\omega_n=-\infty.
\end{aligned}
 \end{equation}
(See \cite{KS} and \cite[chapter 7 and references]{B}. In fact these references treat Hamiltonian cycles. With regard to the case
that $\lim_{n\to\infty}\omega_n=\infty$, it is shown that the limit above holds for Hamiltonian cycles when $k=1$. Since any Hamiltonian cycle can
be cut  open in $n$ possible locations, yielding $n$ Hamiltonian paths, we obtain the result above for any $k$.)

The above result shows in particular that under $P_{n,p_n}$, the Hamiltonian path counting function $N_n^\text{ham}$ converges to the $\delta$-distribution
at $\infty$ if
$p_n$ is as above with $\lim_{n\to\infty}\omega_n=\infty$.
The naive intuition noted in the remark after Theorem \ref{hyper}
 would suggest that the tampering in this case  would be strongly undetectable.
After all, how much can one additional Hamiltonian path be felt in such a situation?
However, we now   demonstrate easily  that whenever $\lim_{n\to\infty} p_n=0$,  the tampering is detectable, while
whenever $p_n\equiv p\in(0,1)$ is constant, the tampering is not strongly  undetectable.
(In fact, it is weakly undetectable, but we will not show that here.)
In light of Proposition \ref{sizebias}, this also shows that when
$p_n\equiv p\in(0,1)$ is constant, the weak law of large numbers does not hold for
$N_n^\text{ham}$, a fact that has been pointed out by Jansen \cite{J}, where a lot of additional results concerning
$N_n^\text{ham}$ can be found.

Label the edges of $G_n$ from 1 to $|e_n|=\frac12n(n-1)$. The random graph   $G(n,p_n)$ with probability measure $P_{n,p_n}$ is constructed
 by considering a collection $\{B_j\}_{j=1}^{|e_n|}$ of IID Bernoulli random variables taking on the values
 1 and 0 with respective probabilities $p_n$ and $1-p_n$, and declaring the $j$-th edge to exist
  if and only if $B_j=1$.  Let $N_n^\text{edges}:\mathcal{E}_n\to \{0,1\cdots, |e_n|\}$
count the number of edges present in an edge configuration. So
under  $P_{n,p_n}$, one has that $N_n^\text{edges}$ is the sum of IID random variables:
$N_n^\text{edges}=\sum_{i=1}^{|e_n|}B_j$.

The expected value of $N_n^\text{edges}$ under the measure $P_{n,p_n}$ is $|e_n|p_n$. Now the tampering involved selecting $n-1$ edges from $e_n$ and demanding that they
exist in the tampered graph. Thus,
the expected value of $N^\text{edges}_n$ under the tampered measure $P_{n,p_n}^\text{ham}$
is  $(|e_n|-(n-1))p_n+(n-1)$.
The increase in  the mean
of $N^\text{edges}_n$ when using the tampered measure instead of the original one is thus
equal to $(1-p_n)(n-1)$.
We denote this change in mean by  $\Delta\text{Exp}_n$.
The variance of  $N^\text{edges}_n$ under the untampered measure is  $|e_n|p_n(1-p_n)$,
and  under  the tampered measure is $(|e_n|-(n-1))p_n(1-p_n)$. Note that these two variances are on the same order since $|e_n|$ is on the order $n^2$.
Let $\text{SD}_n\equiv\sqrt{|e_n|p_n(1-p_n)}$ denote the standard deviation under the untampered measure.
Using the central limit theorem, it is easy to show that if
 $\Delta\text{Exp}_n$ is on a larger order than
  $\text{SD}_n$, then the tampering is
detectable, while if $\Delta\text{Exp}_n$ is on the same order as   $\text{SD}_n$, then the tampering is not strongly undetectable.
In the case that $\lim_{n\to\infty}p_n=0$, we have $\Delta\text{Exp}_n$ on the order $n$ and
$\text{SD}_n$ on the order $o(n)$, while in the case that $p_n=p\in(0,1)$ is constant, we have both
$\Delta\text{Exp}_n$ and
$\text{SD}_n$ on the order $n$.

\medskip

\noindent \it Example 2.\rm\
Consider a random  permutation $\sigma\in S_n$ as a
row of $n$ cards labeled from 1 to $n$ and laid out from left to right in random order. Now tamper with the cards as follows. Select $k_n$ of the cards at random,
remove them from the row, and then replace them in the vacated spaces in increasing order. Let $U_n$
denote the uniform measure on $S_n$, that is, the measure corresponding to a ``random permutation,'' and let
$U_n^{\text{incsubseq},k_n}$ denote the measure on $S_n$ induced from $U_n$ by the above tampering.
Note that by construction, a permutation $\sigma\in S_n$ will have an increasing sequence of length $k_n$ with
$U_n^{\text{incsubseq},k_n}$-probability 1. On the other hand, the celebrated result concerning the length of the longest increasing subsequence
in a random permutation (\cite{LS}, \cite{VK}, \cite{BDJ})  states that the $U_n$ probability of there
being an increasing subsequence of length $cn^\frac12$ goes to 0 as $n\to\infty$, if $c>2$. Thus, one certainly
has $\lim_{n\to\infty}||U_n-U_n^{\text{incsubseq},k_n}||_{\text{TV}}=1$, if $k_n\ge cn^\frac12$, with $c>2$.
The above-mentioned result also states that  the $U_n$-probability  of there being
an increasing subsequence of length $cn^\frac12$ goes to 1 as $n\to\infty$, if $c<2$. From this it follows that
for $k_n\le cn^\frac12$, $c<2$,  the distribution of the  \it number\rm\ of increasing subsequences of length
$k_n$, which we denote by $N_n^{\text{incr},k_n}$, converges to the $\delta$-distribution at $\infty$ as $n\to\infty$. The naive intuition in
the remark after Theorem \ref{hyper} would suggest that
one can tamper on the order $k_n$ without detection, if $k_n\le cn^\frac12$ with $c<2$;
 after all, how much can one additional increasing subsequence be felt in such a situation? However, this turns out to be
false.
In \cite{P06}, it was shown that
$\lim_{n\to\infty}||U_n-U_n^{\text{incsubseq},k_n}||_{\text{TV}}=0$, if $k_n\le n^l$ with $l<\frac25$ and in \cite{P07}
it was shown that $\lim_{n\to\infty}||U_n-U_n^{\text{incsubseq},k_n}||_{\text{TV}}=1$, if $k_n\ge n^l$ with $l>\frac49$.
So in the former case the tampering is strongly undetectable and in the latter case it is detectable.
\medskip

In section 2 we give the proof of Proposition \ref{sizebias} and of part (ii) of Corollary \ref{cheb}.
In section 3 we prove Theorem \ref{hyper}. The proof of parts (a-i) and (b-i) are almost immediate
using the first moment method. The proof of parts (a-ii) (b-ii) and (b-iii) use the second moment method
and involve some quite nontrivial computations, some of which may be interesting in their own right.

\section{Proof of Proposition \ref{2ndmom} and Corollary \ref{cheb}-ii.}
\noindent \it Proof of Proposition \ref{2ndmom}.\rm\
Let $\omega\in \mathcal{E}_n$. Then we have $P_{n,p_n}(\omega \ |O_{n,j})=\frac{1_{\{O_{n,j}\}}(\omega)P_{n,p_n}(\omega)}{P_{n,p_n}(O_{n,j})}$.
Since $N_n(\omega)=\sum_{j=1}^{m_n}1_{O_{n,j}}(\omega)$, and since
 the $O_{n,j}$ have the same $P_{n,p_n}$-probabilities for all $j$,  we have
$E_{n,p_n}N_n=m_nP_{n,p_n}(O_{n,1})$.
Using these  facts along  with the definition of $Q_n$ in \eqref{tamperedmeas}
we have
$$
\begin{aligned}
&Q_n(\omega)=\frac1{m_n}\sum_{j=1}^{m_n} P_{n,p_n}(\omega \ |O_{n,j})=\frac{P_{n,p_n}(\omega)}{m_nP_{n,p_n}(O_{n,1})}\sum_{j=1}^{m_n}
1_{\{O_{n,j}\}}(\omega)=\\
&\frac{N_n}{E_{n,p_n}N_n}P_{n,p_n}(\omega).
\end{aligned}
$$
\hfill $\square$

\noindent \it Proof of Corollary \ref{cheb}-ii.\rm\
Let $Y_n=\frac{N_n}{E_{n,p_n}N_n}$.
Using Proposition \ref{sizebias} along with an  alternative equivalent definition of the total variation distance, we have
$$
||Q_n-P_{n,p_n}||_{\text{TV}}=
\sum_{\omega\in \mathcal{E}_n}(1-\frac{N_n(\omega)}{E_{n,p_n}N_n})^+~P_{n,p_n}(\omega),
$$
where $a^+=a\vee 0$. From this
 it follows that
$\lim_{n\to\infty}||Q_n-P_{n,p_n}||_{\text{TV}}=1$ if and only if
$\lim_{n\to\infty}P_{n,p_n}(Y_n>\epsilon)=0$,
for all $\epsilon>0$.
By the assumption in part (ii) of the corollary,
$E_{n,p_n}Y_n^2\le M$ for some $M$ and all   $n$.
For every $\epsilon>0$,  we have
$$
\begin{aligned}
&1=E_{n,p_n}Y_n\le\epsilon
+E_{n,p_n}Y_n1_{Y_n>\epsilon}\le\epsilon+\big(P_{n,p_n}(Y_n>\epsilon)\big)^\frac12(E_{n,p_n}Y_n^2)^\frac12\le\\
&\epsilon+\big(MP_{n,p_n}(Y_n>\epsilon)\big)^\frac12.
\end{aligned}
$$
From this it is not possible that
$\lim_{n\to\infty}P_{n,p_n}(Y_n>\epsilon)=0$, if $\epsilon<1$.
\hfill $\square$

\section{Proof of Theorem \ref{hyper}}
We begin with the quick proofs of (a-i) and (b-i).

\noindent \it Proof of  (a-i).\rm\
There is a two-to-one correspondence between $H_2^n\times S_n$ and diameter paths in $H_2^n$.
Indeed, for  $\bar x\in H_2^n$ and $\sigma\in S_n$, we begin the diameter path at $\bar x$ and use the permutation
$\sigma$ to determine the order in which we change the components of $\bar x$. (The correspondence is two to one
because the diameter path is not oriented.)
In particular there are $2^{n-1}n!$ diameter paths. The probability that any particular diameter path is contained in
the random hypercube $H_2^n(p_n)$ is $p_n^n$; thus we have
\begin{equation}\label{1hyper}
E_{n,p_n}N_n^{\text{diam}}=2^{n-1}n!p_n^n.
\end{equation}
From this it follows that $\lim_{n\to\infty}E_{n,p_n}N_n^{\text{diam}}=0$,
if $p_n\le \frac \gamma n$, with $\gamma<\frac e2$.
Thus, for such $p_n$, $N_n^\text{diam}$ under $P_{n,p_n}$ converges to the $\delta$-distribution at 0 as $n\to\infty$,
from which it follows that the tampering is detectable.

\noindent \it Proof of (b-i).\rm\ There is a one-to-one correspondence between $S_n$ and diameter paths that start at 0.
The probability that any
particular diameter path is contained in
the random hypercube $H_2^n(p_n)$ is $p_n^n$; thus we have
\begin{equation}\label{1hyper0}
E_{n,p_n}N_n^{\text{diam},0}=n!p_n^n.
\end{equation}
From this it follows that $\lim_{n\to\infty}E_{n,p_n}N_n^{\text{diam},0}=0$,
if $p_n\le \frac \gamma n$, with $\gamma< e$.
As in part (a-i), it then follows that the tampering is detectable.

\medskip

By Corollary \ref{cheb}, to prove (a-ii) it suffices to show that
\begin{equation}\label{chebapplic}
\text{Var}_{n,p_n}(N_n^{\text{diam}})=o\big((E_{n,p_n}N_n^{\text{diam}})^2\big),
\end{equation}
if $p_n$ is as in (a-ii),
and to prove (b-iii) it suffices to show that
\begin{equation}\label{chebapplic0}
\text{Var}_{n,p_n}(N_n^{\text{diam},0})=o\big((E_{n,p_n}N_n^{\text{diam},0})^2\big),
\end{equation}
if $p_n$ is as in (b-iii).

With regard to (b-ii), note that under the untampered  measure, the probability that 0 is an  isolated vertex is
$(1-p_n)^n$. If $p_n$ is as in (b-ii), then this probability stays bounded from 0. On the other hand,
under the tampered measure, the probability that 0 is isolated is 0.
Thus, in the case of (b-ii), the tampering cannot be strongly undetectable. Thus, by Corollary \ref{cheb}, to complete the proof
that the tampering is weakly detectable, it suffices to show that
\begin{equation}\label{genchebapplic0}
\text{Var}_{n,p_n}(N_n^{\text{diam},0})=O\big((E_{n,p_n}N_n^{\text{diam},0})^2\big),
\end{equation}
if $p_n$ is as in (b-ii).

We now give the long and involved proof of \eqref{chebapplic} to prove (a-ii).
After that we will only need a single long paragraph to describe
the changes required to proof \eqref{chebapplic0} and \eqref{genchebapplic0},
which are a bit less involved.

The diameter paths are labeled from 1 to $m_n=2^{n-1}n!$, and we have defined $O_{n,j}$ to be the set of edge
configurations which contain the $j$-th diameter edge.
We relabel for convenience. Let $O_{\bar x,\sigma}$ denote the set of edge configurations
which contain the diameter path corresponding to $(\bar x,\sigma)$ in the above two-to-one correspondence.
Then we have $N_n^{\text{diam}}=\frac12\sum_{\bar x\in H_2^n,\sigma\in S_n}1_{O_{\bar x,\sigma}}$.
Thus
\begin{equation}\label{2hyper}
E_{n,p_n}(N_n^{\text{diam}})^2=\frac14\sum_{\bar x,\bar y\in H_2^n,\sigma,\tau\in S_n}P_n^{p_n}(O_{\bar x,\sigma}\cap
O_{\bar y,\tau}).
\end{equation}
By symmetry considerations, letting id denote the identity permutation and letting $\bar0\in H_2^n$ denote the element with
zeroes in all of its coordinates, we have
\begin{equation}\label{symhyper}
\sum_{\bar x,\bar y\in H_2^n,\sigma,\tau\in S_n}P_n^{p_n}(O_{\bar x,\sigma}\cap
O_{\bar y,\tau})=2^nn!\sum_{\bar x\in H_2^n,\sigma\in S_n}P_n^{p_n}(O_{\bar x,\sigma}\cap
O_{\bar 0,\text{id}}).
\end{equation}
Let
$W_n(\bar x,\sigma)$ denote the number of edges that the diameter path corresponding to $(\bar x, \sigma)$ has in common
with the diameter path corresponding to $(\bar 0,\text{id})$. Then we have
\begin{equation}\label{rvW}
P_n^{p_n}(O_{\bar x,\sigma}\cap
O_{\bar 0,\text{id}})=p_n^{2n-W_n(\bar x,\sigma)}.
\end{equation}
 Letting the generic $E$ denote the expectation
with respect to the uniform measure on $H_2^n\times S_n$,
it then follows from \eqref{1hyper} and  \eqref{2hyper}-\eqref{rvW} that
\begin{equation}\label{key}
E_{n,p_n}(N_n^{\text{diam}})^2=(E_{n,p_n}N_n^{\text{diam}})^2Ep_n^{-W_n}.
\end{equation}
Thus, if we show that
\begin{equation}\label{needhyper}
\lim_{n\to\infty}E(p_n^{-W_n};W_n\ge1)=0,
\end{equation}
then it will follow from \eqref{key} that
\eqref{chebapplic} holds.

We now estimate $P(W_n\ge m)$, for $m\ge1$, where $P$ denotes the probability corresponding to the expectation $E$.
In fact, in the quite involved  estimate
that follows, it will be convenient to assume that $m\ge2$; one can show that the estimate obtained below in
\eqref{finall} also holds for $m=1$.
 The diameter path $(\bar0,\text{id})$ has $n$ edges, which we  label
$e_1,e_2,\cdots, e_n$, with $e_1$ being the edge connecting $\bar 0$ to $(1,0\cdots, 0)$, $e_2$
being the edge connecting $(1,0\cdots, 0)$ to $(1,1,0\cdots, 0)$, etc.
(At the beginning of the paper, $e_n$ was used for the set of edges in $H_2^n$; such use for $e_n$ will not appear again.)
Let $A_{l_1,\cdots, l_m}\subset H_2^n\times S_n$
denote those diameter paths which contain the edges $e_{l_1},\cdots,e_{l_m}$.
Then
\begin{equation}\label{subadd}
P(W_n\ge m)\le \sum_{1\le 1_1<l_2<\cdots <l_m\le n}P(A_{l_1,\cdots, l_m}).
\end{equation}
We now estimate $P(A_{l_1,\cdots, l_m})$.

We first determine  for which $\sigma\in S_n$ one has that $(\bar0,\sigma)\in A_{l_1,\cdots, l_m}$; this result
will be needed  for the general case of determining  which $(\bar x,\sigma)$ belong to $A_{l_1,\cdots, l_m}$.
We will  say that $[j]$ is a sub-permutation of  $\sigma$ if $\sigma$ maps $[j]$ onto itself.
A moment's thought reveals that the edge $e_j$  belongs to the diameter path $(\bar 0,\sigma)$ if and only if both $[j-1]$ and $[j]$ are sub-permutations
for $\sigma$.
Thus,  $(\bar0,\sigma)\in A_{l_1,\cdots, l_m}$ if and only $[l_1-1], [l_1],[l_2-1],[l_2],\cdots, [l_m-1]$ and $[l_m]$ are all sub-permutations
of $\sigma$. The number of permutations $\sigma\in S_n$ for which this holds is easily seen to be
$(l_1-1)!(l_2-1-l_1)!\cdots(l_m-1-l_{m-1})!(n-l_m)!$.
Let $T^n_{m;l_1,\cdots, l_m}\subset S_n$ denote those permutations
for which  $[l_1-1], [l_1],[l_2-1],[l_2],\cdots, [l_m-1]$ and $[l_m]$ are all sub-permutations. So we have
\begin{equation}
|T^n_{m;l_1,\cdots, l_m}|=(l_1-1)!(l_2-1-l_1)!\cdots(l_m-1-l_{m-1})!(n-l_m)!.
\end{equation}

We now consider when
$(\bar x,\sigma)\in A_{l_1,\cdots, l_m}$ for general $\bar x$.
 It is not hard to see that a necessary condition for $(\bar x,\sigma)\in A_{l_1,\cdots, l_m}$ is  that either
$\bar x=(x_1,\cdots, x_n)$ satisfies $x_j=0$, for all $l_1\le j\le l_m$, or $x_j=1$, for
all $l_1\le j\le  l_m$. We will refer to  these two conditions on $\bar x$
by $K_{0;l_1,l_m}$ and $K_{1;l_1,l_m}$.

 If one of these two  conditions on $\bar x$ is satisfied, then in order to have  $(\bar x,\sigma)\in A_{l_1,\cdots, l_m}$,
 the following conditions are required on  $\sigma$. Recall that $\sigma$ gives the order in which the $n$ coordinates of
$\bar x$ are changed so that the diameter path moves from $\bar x$ to $\bar 1-\bar x$. So if $\sigma=(\sigma_1,\cdots, \sigma_n)$,
then the $j$-th edge in the diameter path will involve changing the $\sigma_j$-th coordinate.
Let $B_{0;l_1}(\bar x)$ denote those $j\in\{1,\cdots, l_1-1\}$  for which
$x_j=0$, and let $C_{1;l_m}(\bar x)$ denote those $j\in\{l_m+1,\cdots, n\}$ for
which $x_j=1$ ($B_{0;1}(\bar x),C_{0;m}(\bar x)=\emptyset$).
Let $r_{l_1,l_m}(\bar x)=|B_{0;l_1}(\bar x)|+|C_{1;l_m}(\bar x)|$.
 Then it is not hard to see that  the first
 $r_{l_1,l_m}(\bar x)$ coordinates in $\sigma$ must be reserved for
 $B_{0;l_1}(\bar x)\cup C_{1;l_m}(\bar x)$; that is,
 $\{\sigma_1,\cdots, \sigma_{r_{l_1,l_m}(\bar x)}\}=B_{0;l_1}(\bar x)\cup C_{1;l_m}(\bar x)$.
 Let $x^{1;j}$ denote the vertex in $H_2^n$ whose first $j$ components are 1 and whose remaining components
are 0. Of course, this vertex belongs to the diameter path $(\bar 0,\text{id})$.
If $\sigma$ is as above, then the
$r_{l_1,l_m}(\bar x)$-th vertex of the diameter path $(\bar x,\sigma)$ will be
$x^{1;l_1-1}(\bar x)$ if $\bar x$ satisfies  condition $K_{0;l_1,l_m}$, and will be
$x^{1;l_m}$
if  $\bar x$ satisfies condition  $K_{1;l_1,l_m}$. In the former case, we must then
have $\sigma_{r_{l_1,l_m}(\bar x)+1}=l_1$, and in the latter case, we must then have
$\sigma_{r_{l_1,l_m}(\bar x)+1}=l_m$. In the former case, the
$(r_{l_1,l_m}(\bar x)+1)$-th vertex of the diameter path $(\bar x,\sigma)$ will be
$x^{1;l_1}(\bar x)$ and the $r_{l_1,l_m}(\bar x))$-th edge will be $e_{l_1}$, and in the latter
case,  the
$(r_{l_1,l_m}(\bar x)+1)$-th vertex of the diameter path $(\bar x,\sigma)$ will be
$x^{1;l_m-1}(\bar x)$ and the $r_{l_1,l_m}(\bar x))$-th edge will be $e_{l_m}$.
(Recall that a diameter path has $n+1$ vertices.)

If $\bar x$ satisfies condition  $K_{0;l_1,l_m}$, then the next $l_m-l_1$ coordinates
 of $\sigma$ must involve the numbers $(l_1+1,l_1+2,\cdots, l_m)$, and must move the diameter path $(\bar x, \sigma)$  from
 the vertex $x^{1;l_1}(\bar x)$ to the vertex $x^{1;l_m}(\bar x)$ while passing through
 the edges $e_{l_2},\cdots, e_{l_m}$. Based on our analysis above, for this to happen one requires that
 $(\sigma_{l_1+1}-l_1,\sigma_{l_1+2}-l_1,\cdots, \sigma_{l_m}-l_1)$  belong to
 $T^{l_m-l_1}_{m-2;l_2-l_1,\cdots, l_{m-1}-l_1}\subset S_{l_m-l_1}$.
Similarly, if $\bar x$ satisfies
condition  $K_{1;l_1,l_m}$, then the next $l_m-l_1$ coordinates
 of $\sigma$ must involve the numbers $(l_1+1,l_1+2,\cdots, l_m)$, and must move the diameter path $(\bar x, \sigma)$  from
 the vertex $x^{1;l_m}(\bar x)$ to the vertex $x^{1;l_1}(\bar x)$ while passing through
 the edges $e_{l_{m-1}},\cdots, e_{l_1}$. Inverting the direction of   our analysis above, for this to happen one requires that
 $(\sigma_{l_1+1}-l_1,\sigma_{l_1+2}-l_1,\cdots, \sigma_{l_m}-l_1)$  belong to
 $T^{l_m-l_1}_{m-2;l_m-l_{m-1},\cdots, l_m-l_2}\subset S_{l_m-l_1}$.
Then
 finally, the last $n-r_{l_1,l_m}(\bar x)-1-(l_m-l_1)$ coordinates of $\sigma$ can be chosen arbitrarily from
 the remaining numbers.
Putting the above all together, we obtain
\begin{equation}\label{bigcalc}
\begin{aligned}
&P(A_{l_1,\cdots, l_m})=\\
&\frac1{2^nn!}\sum_{c=0}^{n-l_m}\sum_{b=0}^{l_1-1}\binom{n-l_m}c\binom{l_1-1}b(b+c)!(n-b-c-1-(l_m-l_1))!\times\\
&(|T^{l_m-l_1}_{m-2;l_2-l_1,\cdots, l_{m-1}-l_1}|+|T^{l_m-l_1}_{m-2;l_m-l_{m-1},\cdots, l_m-l_2}|).
\end{aligned}
\end{equation}
(Given that $\bar x$ satisfies condition $K_{0;l_1,l_m}$ or condition $K_{1;l_1,l_m}$,
there are  $\binom{n-l_m}c\binom{l_1-1}b$ ways to choose $\bar x$ so that
$b=|B_{0;l_1}(\bar x)|$ and $c=|C_{1;l_m}(\bar x)|$. And given this, there are
$(b+c)!(n-b-c-1-(l_m-l_1))!(|T^{l_m-l_1}_{m-2;l_2-l_1,\cdots, l_{m-1}-l_1}|$ ways to choose $\sigma$ if condition
$K_{0;l_1,l_m}$ was satisfied, and $(b+c)!(n-b-c-1-(l_m-l_1))!|T^{l_m-l_1}_{m-2;l_m-l_{m-1},\cdots, l_m-l_2}|$
ways to choose $\sigma$ if condition $K_{1;l_1,l_m}$ was satisfied.)

We have
\begin{equation}\label{bcsum}
\begin{aligned}
&\sum_{c=0}^{n-l_m}\sum_{b=0}^{l_1-1}\binom{n-l_m}c\binom{l_1-1}b(b+c)!(n-b-c-1-(l_m-l_1))!=\\
&\sum_{c=0}^{n-l_m}\sum_{b=0}^{l_1-1}\frac{\binom{n-l_m}c\binom{l_1-1}b}{\binom{n-1-l_m+l_1}{b+c}}(n-1-(l_m-l_1))!\le n^2(n-1-(l_m-l_1))!,
\end{aligned}
\end{equation}
where the last inequality follows from the fact that the fraction in the sum above is always less than 1.
After completing the current proof, we will prove the following proposition.
\begin{proposition}\label{keyhyper}
For every $\delta>0$, there exist a $c_\delta>0$ and an $r_\delta\ge0$ such that
$|T^n_{m;l_1,\cdots, l_m}|\equiv(l_1-1)!(l_2-1-l_1)!\cdots(l_m-1-l_{m-1})!(n-l_m)!$ satisfies
$$
\sum_{1\le l_1<l_2<\cdots<l_m\le n}|T^n_{m;l_1,\cdots, l_m}|\le c_\delta m^{r_\delta}(1+\delta)^m(n-m)!,\  1\le m\le n<\infty.
$$
\end{proposition}
From Proposition \ref{keyhyper}, it follows that for any $\delta>0$, there exists a $c_\delta>0$ and an $r_\delta\ge0$ such that
\begin{equation}\label{suml's}
\begin{aligned}
&\sum_{l_1<l_2<\cdots l_{m-1}<l_m}(|T^{l_m-l_1}_{m-2;l_2-l_1,\cdots, l_{m-1}-l_1}|+|T^{l_m-l_1}_{m-2;l_m-l_{m-1},\cdots, l_m-l_2}|)\le\\
&2c_\delta m^{r_\delta}(1+\delta)^m(l_m-l_1-m+2)!.
\end{aligned}
\end{equation}
(Note that in the sum above,  the last subscript, $l_{m-1}-l_1$  in $T^{l_m-l_1}_{m-2;l_2-l_1,\cdots, l_{m-1}-l_1}$
and  $l_m-l_2$ in $T^{l_m-l_1}_{m-2;l_m-l_{m-1},\cdots, l_m-l_2}$, is strictly less than the superscript $l_m-l_1$, whereas
in the sum in Proposition \ref{keyhyper} the last subscript, $l_m$ in  $T^n_{m;l_1,\cdots, l_m}$, can attain the value $n$ of the superscript; however,
this is no problem since the inequality goes in the right direction.)
Now \eqref{bigcalc}, \eqref{bcsum} and \eqref{suml's} give
\begin{equation}
\begin{aligned}
&\sum_{l_1<l_2<\cdots<l_{m-1}<l_m}P(A_{l_1,\cdots, l_m})\le\\
& \frac{2c_\delta m^{r_\delta}n^2}{2^n}(1+\delta)^m
\frac{(n-1-(l_m-l_1))!(l_m-l_1-m+2)!}{n!}.
\end{aligned}
\end{equation}
Now summing over $l_1$ and $l_m$, and denoting $k=l_m-l_1+1$, we have
\begin{equation}\label{almostthere}
\sum_{1\le l_1<l_2<\cdots<l_{m-1}<l_m\le n}P(A_{l_1,\cdots, l_m})\le
\frac{2c_\delta m^{r_\delta} n^3}{2^n}(1+\delta)^m\sum_{k=m}^n\frac1{\binom nk}\frac{(k-m+1)!}{k!}.
\end{equation}
Let $\rho(k)\equiv\frac1{\binom nk}\frac{(k-m+1)!}{k!}=\frac{(k-m+1)!}{n(n-1)\cdots (n-k+1)}$, $m\le k\le n$,
and let $h(k)=\frac{\rho(k+1)}{\rho(k)}$. It is easy to check that $h$ is increasing, which implies that $\rho$ is convex.
Thus, $\rho$ attains its maximum at an endpoint. We conclude that
the maximum of $\rho(k)$ is $\rho(n)=\frac{(n-m+1)!}{n!}$. Using Stirling's formula, it is easy to check that
there exists a $K$ such that $\frac{(n-m+1)!}{n!}\le K(\frac en)^{m-1}$. Using these facts in
\eqref{almostthere}, we obtain
\begin{equation}\label{almostfinal}
\sum_{1\le l_1<l_2<\cdots<l_{m-1}<l_m\le n}P(A_{l_1,\cdots, l_m})\le\frac{2Kc_\delta m^{r_\delta} n^5}{2^ne}\big(\frac{(1+\delta)e}n\big)^m.
\end{equation}

Using \eqref{almostfinal} in \eqref{subadd} now gives
\begin{equation}\label{finall}
P(W_n\ge m)\le \frac{2Kc_\delta m^{r_\delta} n^5}{2^ne}\big(\frac{(1+\delta)e}n\big)^m.
\end{equation}
Thus, if $p_n=\frac\gamma n$, then from \eqref{finall} we have
\begin{equation}\label{reallyfinal}
E(p_n^{-W_n}; W\ge 1)\le  \sum_{m=1}^n (\frac n\gamma)^mP(W_n\ge m)\le \frac{2Kc_\delta n^{r_\delta+5}}{2^ne}\sum_{m=1}^n(\frac{(1+\delta)e}{\gamma})^m.
\end{equation}
We may choose $\delta>0$ as small as we like in \eqref{reallyfinal}.
For  $\gamma>\frac e2$,  choose $\delta$ so that $\frac{(1+\delta)e}\gamma<2$.
Then it follows from \eqref{reallyfinal} that
\eqref{needhyper} holds for $p_n=\frac\gamma n$ with $\gamma>\frac e2$.
\hfill $\square$

We now return to prove Proposition \ref{keyhyper}.

\noindent \bf Proof of Proposition \ref{keyhyper}.\rm\ Define
$$
C^{(0)}_j=\sum_{i=0}^j\frac1{\binom ji},\ j\ge0,
$$
and then define by induction the iterates
$$
C_j^{(m)}=\sum_{i=0}^j\frac{C^{(m-1)}_i}{\binom ji},\ j\ge0,m\ge1.
$$
We have
\begin{equation}\label{0}
\sum_{1\le l_1<l_2}(l_1-1)!(l_2-1-l_1)!=(l_2-2)!\sum_{l_1=1}^{l_2-1}\frac1{\binom{l_2-2}{l_1-1}}=(l_2-2)!C^{(0)}_{l_2-2},
\end{equation}
and then using \eqref{0},
\begin{equation*}
\begin{aligned}
&\sum_{1\le l_1<l_2<l_3}(l_1-1)!(l_2-1-l_1)!(l_3-1-l_2)!=\\
&\sum_{2\le l_2<l_3}(l_2-2)!C^{(0)}_{l_2-2}(l_3-1-l_2)!=
(l_3-3)!\sum_{l_2=2}^{l_3-1}\frac{C^{(0)}_{l_2-2}}{\binom{l_3-3}{l_2-2}}=(l_3-3)!C^{(1)}_{l_3-3}.
\end{aligned}
\end{equation*}
Continuing in this vein, we obtain
\begin{equation*}
\sum_{1\le l_1<l_2<\cdots<l_m}(l_1-1)!(l_2-1-l_1)!\cdots(l_m-1-l_{m-1})!=
(l_m-m)!C_{l_m-m}^{(m-2)},
\end{equation*}
and
\begin{equation}\label{keyprop}
\sum_{1\le l_1<l_2<\cdots<l_m<n}(l_1-1)!(l_2-1-l_1)!\cdots(l_m-1-l_{m-1})!(n-l_m)!=(n-m)!C_{n-m}^{(m-1)}.
\end{equation}
In light of \eqref{keyprop}, to complete the proof of Proposition \ref{keyhyper}, it suffices to show that
for every $\delta>0$, there exist a  $c_\delta>0$ and an $r_\delta\ge0$ such that
\begin{equation}\label{needprop}
\sup_{n\ge1} C_n^{(k)}\le c_\delta k^{r_\delta}(1+\delta)^k,\ k\ge1.
\end{equation}

Let $n_0\ge1$, and for $n>n_0$ write
\begin{equation}\label{recurs}
C_n^{(k)}=\sum_{i=0}^{n_0}\frac{C_i^{(k-1)}}{\binom n i}+\sum_{i=n_0+1}^n\frac{C_i^{(k-1)}}{\binom n i},\ k\ge1, n>n_0.
\end{equation}
We need the following lemma whose proof we defer until the completion of the proof of the proposition.
\begin{lemma}\label{lemma}
For each $n$ there exists  a constant $c_n$ such that
\begin{equation}\label{lemmamatrix}
C^{(k)}_n\le c_nk^n,\  k\ge1.
\end{equation}
\end{lemma}
From \eqref{recurs} and \eqref{lemmamatrix}, it follows that for each $n_0$ there exists a constant $\gamma_{n_0}$ such that
\begin{equation}\label{recursimprove}
C_n^{(k)}\le \gamma_{n_0}(k-1)^{n_0}+\sum_{i=n_0+1}^n\frac{C_i^{(k-1)}}{\binom n i},\ k\ge1, n>n_0.
\end{equation}
Let
\begin{equation}\label{defdn}
d_{n_0}\equiv\sup_{n\ge n_0+1}\sum_{i=n_0+1}^n\frac1{\binom ni}.
\end{equation}
It is easy to see that
\begin{equation}\label{easylim}
\lim_{n_0\to\infty}d_{n_0}=1.
\end{equation}
Letting
$$
A_{n_0}^{(k)}\equiv\sup_{n>n_0}C_n^{(k)},
$$
we have from \eqref{recursimprove} and \eqref{defdn} that
\begin{equation}\label{keyrecurs}
A_{n_0}^{(k)}\le \gamma_{n_0}(k-1)^{n_0}+d_{n_0}A_{n_0}^{(k-1)},\ k\ge1.
\end{equation}
It is not hard to show that
$$
\sup_{n\ge0}C_n^{(0)}=\sup_{n\ge0}\sum_{i=0}^n\frac1{\binom ni}=\frac83;
$$
however, all we need for our purposes is that this quantity is bounded, and this is very easy  to see.
Thus, we have
\begin{equation}\label{keyrecurs0}
A_{n_0}^{(0)}\le \frac83.
\end{equation}
It is easy to show that if $\{x_j\}_{j=0}^k$ satisfies the recursive inequalities
$x_0\le \frac83$ and $x_j\le C+d_{n_0}x_{j-1}$, for $1\le j\le k$, then
$x_k\le C(1+d_{n_0}+\cdots+d_{n_0}^{k-1}+\frac83d_{n_0}^k)=C\big(\frac{d_{n_0}^k-1}{d_{n_0}-1}+\frac83d_{n_0}^k\big)$.
Applying this with $C=\gamma_{n_0}(k-1)^{n_0}$, it follows from \eqref{keyrecurs} and \eqref{keyrecurs0} that
\begin{equation}\label{finalhyper}
\sup_{n>n_0}C_n^{(k)}=A_{n_0}^{(k)}\le \gamma_{n_0}(k-1)^{n_0}\big(\frac{d_{n_0}^k-1}{d_{n_0}-1}+\frac83d_{n_0}^k\big), \ k\ge1.
\end{equation}
By \eqref{easylim}, for any $\delta>0$, there exists an $n_0$ such that
$d_{n_0}\le 1+\delta$. Using this with \eqref{finalhyper}, and using \eqref{lemmamatrix} with $n\le n_0$,
one concludes that  \eqref{needprop} holds. This completes the proof of the
proposition. \hfill $\square$

We now return to prove Lemma \ref{lemma}.

\noindent \bf Proof of Lemma \ref{lemma}.\rm\
Fix $n\ge1$. Let $B$ denote the $n\times n$ matrix with entries $b_{ij},~ 1\le i,j\le n$, given
by $b_{ij}=\frac1{\binom ij}$, for $i\ge j$, and $b_{ij}=0$, for $j>i$.
Let $v^0$ denote the $n$-vector with entries $v^0_j$, $1\le j\le n$, given by $v^0_j=C^{(0)}_j=\sum_{l=0}^j\frac1{\binom jl}$.
Then from the recursive definition of the $\{C_m^{(k)}\}_{m,k=0}^\infty$, it follows that
\begin{equation}\label{matrix}
C_n^{(k)}=(B^kv^0)_n,\ k\ge1,
\end{equation}
 where $(B^kv^0)_n$ denotes the $n$-th coordinate of the $n$-vector $B^kv^0$.
Since $B$ is lower triangular with all ones on the diagonal, it follows that there exist vectors $v^1,\cdots, v^{n-1}$ such that
$Bv^0=v^0+v^1, Bv^1=v^1+v^2, \cdots, Bv^{n-2}=v^{n-2}+v^{n-1}$ and $Bv^{n-1}=v^{n-1}$.
From this, it follows that
\begin{equation}\label{almostfinalprop}
B^kv^0=\sum_{l=0}^{k\wedge n} \binom klv^l.
\end{equation}
Thus, from \eqref{matrix} and \eqref{almostfinalprop}, we obtain
\begin{equation}\label{finalll}
C_n^{(k)}=\sum_{l=0}^{k\wedge n} \binom klv^l_n.
\end{equation}
where $v^l_n$ is the $n$-th coordinate of $v^l$.
The lemma follows immediately from \eqref{finalll}.\hfill $\square$

We have now completed the proof of \eqref{chebapplic}, and thus the proof of (a-ii).
To complete the proof of (b-ii) and (b-iii) we need to prove \eqref{chebapplic0} and \eqref{genchebapplic0}.
In fact all the work has been done in the above proof. The proof up to \eqref{subadd} is the same as before,
except that now we work with the space $S_n$ instead of with $H_2^n\times S_n$. In particular then,
we now have $W_n=W_n(\sigma)$, and it denotes the number of edges that the diameter path starting from 0
and corresponding to $\sigma$ has in common with the diameter path starting from 0 and corresponding to $\text{id}$.
Similarly, $A_{l_1,\cdots, l_m}\subset S_n$ denotes the number of diameter paths starting from 0 which contain the
edges $e_{l_1},\cdots, e_{l_m}$.
From the paragraph after \eqref{subadd}, it follows that
$A_{l_1,\cdots, l_m}=T^n_{m;l_1,\cdots, l_m}$ and that
\begin{equation}\label{AT}
P(A_{l_1,\cdots, l_m})=
\frac{|T^n_{m;l_1,\cdots, l_m}|}{n!}.
\end{equation}
Using \eqref{AT} with \eqref{subadd} and Proposition \ref{keyhyper}, it follows that
\begin{equation}\label{Wb}
P(W_n\ge m)\le c_\delta m^{r_\delta}(1+\delta)^m\frac{(n-m)!}{n!}.
\end{equation}
As noted above, there exists a $K$ such that $\frac{(n-m)!}{n!}\le K(\frac en)^m$.
Thus, we have
\begin{equation}\label{final}
E(p_n^{-W};W\ge 1)\le Kc_\delta\sum_{m=1}^n(p_n)^{-m}m^{r_\delta}(1+\delta)^m(\frac en)^m.
\end{equation}
From \eqref{final}, it follows that as $n\to\infty$,
$E(p_n^{-W};W\ge 1)$ converges to 0 if $p_n$ is as in  (b-iii),  and remains bounded if $p_n$ is as in (b-ii). Thus, it follows from
\eqref{key}
that \eqref{chebapplic0} holds if $p_n$ is as in (b-iii) and that  \eqref{genchebapplic0} holds if $p_n$ is as in (b-ii).
\hfill $\square$


\begin{thebibliography}{99}

\bibitem{BDJ}
 Baik, J., Deift, P. and Johansson, K.
\emph{On the distribution of the length of the longest increasing subsequence of
random permutations},   J. Amer. Math. Soc.  \textbf{12}  (1999), 1119-1178.


\bibitem{B} Bollab\'as, B.,  \emph{Modern Graph Theory},  Graduate Texts in Mathematics, 184, Springer-Verlag (1998).




\bibitem{J} Janson, S.,  \emph{The numbers of spanning trees, Hamilton cycles and perfect matchings in a random graph},
Combin. Probab. Comput. 3 (1994),  97-126.






\bibitem{KS}
 Koml\'os, J. and  Szemer\'edi, E.,  \emph{Limit distribution for the existence of Hamiltonian cycles in a random graph},
 Discrete Math. 43 (1983),  55-63.





\bibitem{LS}  Logan, B. F. and Shepp, L. A.
\emph{A variational problem for random Young tableaux},
  Advances in Math.  \textbf{26}  (1977), 206-222.



\bibitem{P06}
Pinsky, R. G. , \emph{Law of large numbers for increasing subsequences of random permutations}, Random Structures Algorithms \textbf{29} (2006),  277-295.

\bibitem{P07}
 Pinsky, R. G., \emph{When the law of large numbers fails for increasing subsequences of random permutations}, Ann. Probab. \textbf{35} (2007),  758-772.



\bibitem{VK}
Vershik, A. M. and Kerov, S. V.
\emph{Asymptotic behavior of the maximum and generic dimensions of irreducible
representations of the symmetric group},
 Functional Anal. Appl.  \textbf{19}  (1985),  21-31.



 \end{thebibliography}
\end{document}